\input amstex
\input xy
\input epsf
\xyoption{all}
\documentstyle{amsppt}
\document
\magnification=1200
\NoBlackBoxes
\nologo
\hoffset1.5cm
\voffset2cm
\def\D{\Cal{D}}
\def\C{\bold{C}}
\def\G{\bold{G}}
\def\Q{\bold{Q}}
\def\P{\bold{P}}
\def\Z{\bold{Z}}

\def\g{\gamma}
\pageheight {16.cm}


\medskip

\centerline{\bf MODULAR FORMS OF REAL WEIGHTS}
\smallskip
\centerline{\bf AND GENERALIZED DEDEKIND SYMBOLS} 
\bigskip

\centerline{\bf Yuri I. Manin}

\medskip

\centerline{\it Max--Planck--Institut f\"ur Mathematik, Bonn, Germany}

\bigskip

{\it ABSTRACT.} In a previous paper,  I have defined non--commutative generalised Dedekind symbols
for classical  $PSL(2,\Z )$--cusp forms  using iterated period polynomials. Here I generalise this
construction to forms of real weights using their iterated period functions 
introduced and studied in a recent  article  by R.~Bruggeman and Y.~Choie.
\bigskip

\centerline{\bf 1. Introduction: generalized Dedekind symbols}

\medskip

 The classical Dedekind symbol encodes an essential
part of  modular properties of the Dedekind eta--function, and appears
in many contexts seemingly unrelated to modular forms  (cf. [KiMel], [Mel]).
Fukuhara in [Fu1], [Fu2], and others ([Ap], [ChZ], gave an abstract definition
of generalized  Dedekind symbols with values in an arbitrary commutative group
and produced such symbols from period polynomials of $PSL (2,\Z )$--modular forms of 
any even weight.
\smallskip
In the note [Ma5], I have  given an abstract definition of generalised  Dedekind symbols  
 for the full modular group $PSL (2,\Z )$  taking values in arbitrary  non--necessarily commutative
group and constructed such symbols from iterated versions of period integrals of modular
forms of integral weights considered earlier in [Ma3], [Ma4].
\smallskip
 In this article, I extend these constructions to cusp forms of real weights, 
 studied in particular  in [Kn], [KnMa], [BrChDi]. The essential ingredient here
 is furnished by the introduction of iterated versions of their period integrals following  [BrCh] .
 \smallskip

{\bf 1.1. Dedekind symbols.} The Dedekind eta function is a holomorphic function
of the complex variable $z$ with positive imaginary part given by
$$
\eta (z) = e^{\pi iz/24} \prod_{k=1}^{\infty} (1-e^{2\pi i kz}) .
$$
It is a cusp form of weight $1/2$, and from $PSL (2,\Z )$--invariance
of $\eta (z)^{24} dz^6$ it follows that for any 
fractional linear transformation $\gamma \in PSL (2,\Z )$,
$$
\gamma z := \frac{az+b}{cz+d}
$$
with $c> 0$ we can obtain a rational number
$$
s(a,c):=- \frac{1}{\pi i} [\roman{log}\,\eta (\gamma z) - \roman{log} \eta (z) - \frac{1}{2}
\roman{log} ((cz+d)/i)] + \frac{a+d}{12c}
$$
called the (classical) Dedekind symbol.

\smallskip

It satisfies the reciprocity relation (which can be easily extended to all $c\in \Z$)
$$
s(a,c)+ s(c,a) =a/c+c/a +1/ac -3\, \roman{sgn} (a) .
$$

{\bf 1.2. Generalised Dedekind symbols with values in an abelian group.} Using a slightly different
normalisation and terminology of [Fu1], one can define the generalised Dedekind symbol $d(p,q)$ as a function 
 $d:\,W\to \bold{G}$ where $W$ is the set of pairs of co--prime integers $(p,q)$, and $\bold{G}$ an abelian group. It can be uniquely reconstructed from the
functional equations
$$
d(p,q)=d(p,q+p), \quad d(p,-q)= -d(p,q),
\eqno(1.1)
$$
$$
d(p,q)-d(q,-p) = \frac{p^2+q^2-3pq+1}{12pq}.
\eqno(1.2)
$$
Studying other $PSL(2,\Z )$--modular forms in place of $\eta$, one arrives to {\it generalised Dedekind
symbols}, satisfying  similar
functional equations, in which the right hand side of (1.2) is replaced
by a different {\it reciprocity function}, which in turn satisfies  simpler functional
equations and which uniquely defines the respective Dedekind symbol: see [Fu1], [Fu2] and 1.3.2 below.
\smallskip
In particular, let $F (z)$ be a cusp form of even integral weight $k+2$ for $\Gamma := PSL(2,\bold{Z})$. Its period polynomial
is the following function of $t\in \bold{C}$:
$$
P_F (t):= \int_0^{i\infty} F (z) (z-t)^kdz
\eqno(1.3)
$$

Fukuhara has shown that (slightly normalized) values of period polynomials at rational points
form a reciprocity function
$$
f_F (p, q):= p^k P_F (q/p) = \int_0^{i\infty} F (z) (pz-q)^k dz
\eqno(1.4)
$$
whose respective  Dedekind symbol is
$$
d_F (p,q) =  \int_{p/q}^{i\infty} F (z) (pz-q)^k dz
\eqno(1.5)
$$
\medskip
{\bf 1.3.  Non--commutative generalised Dedekind symbols.} In [Ma5], I introduced {\it non--commutative generalised Dedekind symbols} with values in
a non--necessarily abelian group $\G$ by the  following definition.

\medskip

{\bf 1.3.1. Definition.} {\it A $\bold{G}$--valued reciprocity function  is a map $f: W\to \G$ satisfying the  
following conditions}
$$
f(p,-q)=f(-p,q).
\eqno(1.6)
$$
$$
f(p,q)f(-q,p)=1_{\G}
\eqno(1.7)
$$
$$
f(p,p+q)f(p+q,q)=f(p,q).
\eqno(1.8)
$$
\medskip

Applying (1.8) to $p=1, q=0$, we get $f(1,1)=1_{\bold{G}}$ where $1_{\bold{G}}$ is the identity.
From (1.7)
we then get $f(-1,1)=1_{\bold{G}}.$
Moreover,  $f(-p,-q)=f(p,q)$
so that $f(p,q)$ depends only on $q/p$ which obviously can now be an arbitrary point in $\bold{P}^1(\bold{Q})$
including $\infty$. (Of course, $i\infty$ in integrals like (1.3)--(1.5) coincides with $\infty$ of the real projective line).

\smallskip

The function (1.4) taking values in the additive group of complex
numbers  satisfies equations (1.6)--(1.8) (written additively). 

\medskip

{\bf 1.3.2. Definition.} {\it Let $f$ be a $\bold{G}$--valued reciprocity function.

\smallskip
 A generalized  $\G$--valued Dedekind symbol $\Cal{D}$ with reciprocity
 function $f$ is a map
$$
\Cal{D}:\,W\to \G:\ (p,q)\mapsto \Cal{D}(p,q)
$$
satisfying the following conditions (1.9)--(1.11): 
$$
  \D(p,q)=\D(p, q+p),
\eqno(1.9)
$$
$$
 \D(p,-q) =\D(-p,q),
\eqno(1.10)
$$
so that $\D(-p,-q)= \D(p,q)$.
Finally, 
$$
\D(p,q)\D(q,-p)^{-1}=f(p,q).
\eqno(1.11)
$$
}

Clearly, knowing $\Cal{D}$ one can uniquely reconstruct its reciprocity function $f$.
Conversely, any reciprocity function uniquely defines the respective generalised 
Dedekind symbol ([Ma5], Theorem 1.8).

\medskip

In [Ma5] I constructed such reciprocity functions using iterated integrals of cusp forms
of integral weights. In the main part of this note I will (partly) generalise this construction to cusp forms
of real weights.

\smallskip

Period polynomials of cusp forms of integer weights appear in many interesting contexts.
Their coefficients are values of certain $L$--functions in integral points
of the critical strip ([Ma1], [Ma2] and many other works); they can be used in order to produce
``local zeta--factors'' in the mythical algebraic geometry of  characteristic 1 ([Ma6]); they describe relations
between certain inner derivations of a free Lie algebra ([Po], [Hai], [HaiMo], [BaSch]), essentially because
iterated period polynomials define representations of unipotent completion
of basic fundamental modular groupoids. 

\smallskip

Iterated period polynomials of real weights can be compared to various other constructions
where interpolation from integer values to real values occurs, e.~g. Deligne's
theory of  `` symmetric groups  $S_w$, $w\in \bold{R}$" using a categorification. It would be very interesting
to find  similar categorifying constructions also in the case of modular forms of real weights. One can expect perhaps
appearance of ``modular spaces $\overline{M}_{1,w}, w\in  \bold{R}.$'' Notice that certain
$p$--adic interpolations appeared already long time ago in the theory
of $p$--adic $L$--functions.

\medskip

{\it Acknowledgements.} This note was strongly motivated and inspired by the recent preprint [BrCh] due to
R.~Bruggeman and  Y.~Choie. R.~Bruggeman kindly answered my questions,
and clarified for me many issues regarding modular forms of non--integer weights.
Together with Y.~Choie, he  carefully read a preliminary version of this note.
I am very grateful to them.

\bigskip

\centerline{\bf 2. Modular forms of real weight and their period integrals}

\medskip

In this section I fix notation and give a brief survey of relevant definitions and results from [Kn], [KnMa],
and [BrCh]. I adopt conventions of [BrCh], where modular forms of real weights are holomorphic
functions on the upper complex half-plane, whereas their period integrals, analogs of (1.3),
are holomorphic functions on the lower half--plane.
\medskip
{\bf 2.1. Growth conditions for holomorphic functions in upper/lower complex half--planes.}
Let $\P^1(\C)$ be the set of $\C$--points of the projective line endowed with a fixed projective coordinate
$z$. This coordinate  identifies the complex plane $\C$ with the maximal subset of $\P^1(\C)$ where $z$ is holomorphic.
\smallskip
We put 
$$
H^+ := \{z\in \C\,|\, \roman{Im}\,z>0 \} ,\quad  H^- := \{t\in \C\,|\, \roman{Im}\,t<0 \}.
$$
As in [BrCh], we identify holomorphic functions on $H^+$ and  $H^-$ using {\it antiholomorphic} involution:
$F(z)\mapsto \overline{F(\overline{z})}$. In the future holomorphic functions on $H^-$ will
often be written using coordinate $t=\overline{z}$. On the other hand, the standard hyperbolic
metric of curvature $-1$ on $H^+\cup\, H^-$, $ds^2=|dz|^2/(\roman{Im}\,z)^2$ looks identically in both
coordinates.

\smallskip

Cusps $\P^1(\Q )$ are rational points  on the common boundary of $H^+$ and $H^-$ (including infinite point).
Denote by $\Cal{P}^+$, resp.  $\Cal{P}^-$ the space of functions $F(z)$ holomorphic in $H^+$, resp. $H^-$ 
satisfying for some constants $K, A>0$ and all $z\in H^{\pm}$ inequality
$$
|F(z)|< K(|z|^A+ |\roman{Im} z|^{-A} ) .
\eqno(2.1)
$$
This is called the {\it polynomial growth} condition.

\smallskip

Cusp forms and their iterated period functions with which we will be working actually satisfy stronger growth conditions
near the boundary: see 2.4 below.
\medskip

{\bf 2.2. Actions of modular group.} The standard left action of $PSL(2,\Z )$ upon $\P^1(\C )$
by linear fractional transformations of $z$ 
$$
z\mapsto \gamma (z)= \frac{az+b}{cz+d}
$$
defines the right action upon holomorphic
functions  in $H^{\pm}$. In the theory of automorphic forms of even integral weight 
this natural right action is considered first upon tensor powers of 1--forms $F(z) (dz)^{k/2}$ and then
transported back to holomorphic functions by dividing the result by $(dz)^{k/2}$.
Equivalently, the last action on functions can be defined using integral powers of $j(\gamma , z):=cz+d$:
$$
(F|_k^1\gamma )(z)= j(\gamma ,z)^{-k} F(\gamma z)
\eqno(2.2)
$$
and similarly in $t$--coordinate.

\smallskip

In the theory of automorphic forms of general real weight  $k$, the relevant generalisation requires two additional
conventions. First, we define $(cz+d)^k$ using the following choice of arguments:

$$
\roman{arg}\,(cz+d) \in (-\pi ,\pi ]\  \roman{for}\  z\in H^+,
$$
$$
\roman{arg}\,(ct+d) \in [-\pi ,\pi )\  \roman{for}\  t\in H^-.
$$
Second, multiplication by $(cz+d)^k$ is completed by an additional complex factor
depending on $\gamma$.
\medskip

{\bf 2.2.1. Definition.} {\it  A unitary multiplier system $v$ of weight $k\in \bold{R}$ (for the group $SL(2,\Z)$) is a
map $v:\, SL(2,\Z) \to \C$, $|v(\gamma )|=1$, satisfying the following conditions. Put 
$$
j_{v,k}(\gamma , z):= v(\gamma ) (cz+d )^k.
$$
Then we have
$$
j_{v,k}(\gamma\delta ,z)= j_{v,k}(\gamma, \delta z)\cdot j_{v,k}(\delta ,z) 
\eqno(2.3)
$$
and
$$
j_{v,k}(-\gamma , z)= j_{v,k}(\gamma , z) .
\eqno(2.4)
$$
}

\smallskip
Identities (2.3) and (2.4) imply that the formula
$$
F(z) \mapsto v(\gamma ) j(\gamma ,z)^{k} F(\gamma z)=j_{v,k}(\gamma ,z) F(\gamma z)
\eqno(2.5)
$$
defines a right action of $SL(2,\Z)/(\pm \roman{id})=PSL(2,\Z)$ upon functions holomorphic in $H^+$.
Functions invariant
with respect to this action and having exponential decay at cusps (in terms of geodesic distance,
cf. [KnMa] and {BrChDi]) are called {\it cusp forms} for the full modular group
of weight $k$ with multiplier system $v$.

\smallskip

For such a form $F(z)$, one can define its {\it period function} $P_F(t)$  by the formula similar to
(1.3). Generally, it is defined only on $H^-$ and satisfies the polynomial growth condition
near the boundary.

\smallskip
Moreover, behaviour of period functions of modular forms with respect to modular transformations involves
the action
$$
(P|_{-k}^v\gamma )(t) :=  v(\gamma ) ^{-1}j(\gamma ,t)^{k} P(\gamma t)=j_{v,k}(\gamma ,t)P(\gamma t)
\eqno(2.6)
$$
which is  the right action of $PSL(2,\Z)$ upon functions holomorphic in $H^-$.

\medskip

{\bf 2.3. Modular forms.} Let $F\in \Cal{P}^+$ be a holomorphic function of polynomial growth in $H^+$ (see 2.1)
satisfying the $SL(2,\Z)$--invariance condition:
$$
(F |_{k+2}^v\gamma )(z)  = F (z)\ \roman{for\ all}\ \gamma \in  SL(2,\Z).
\eqno(2.7)
$$
It is called {\it a modular form of weight $k+2$ and multiplier system $v$}. Such a modular form
is called  {\it a cusp form} if in addition
its Fourier series at all cusps contain only positive powers of the relevant
exponential function (cf. [KnMa]).

The space of all such forms is denoted $C^0(\Gamma ,k+2,v)$. It can be non--trivial
only if $k>0$.

\smallskip

{\bf  2.4. Period integrals.} For a cusp  form  $F\in C^0(\Gamma ,k,v)$ and points $a,b$ in $H^+\cup \{ cusps\}$ we 
put $\omega_F(z;t):=F (z)(z-t)^kdz$ and define its  integral
as a function of $t\in H^-$:
$$
I_a^b (\omega_F;t) :=\int _a^b \omega_F (z;t) .
$$
If  $a$ and/or $b$ is a cusp, then the integration path near it must follow a segment of geodesic
connecting $a$ and $b$. We may and will assume that in our (iterated) integrals the 
integration path is always the segment of geodesic connecting limits of integration.

\smallskip

More generally, for a finite sequence of cusp forms $F_j\in  C^0(\Gamma ,k_j+2,v_j)$ and $\omega_j(z) 
= \omega_j(z;t) := F_j(z)(z-t)^{k_j}dz$,
$j=1,\dots ,n$, where $t$ is considered as parameter, we put
$$
I_a^b (\omega_1,\dots ,\omega_n; t) :=\int _a^b \omega_1(z_1)   \int _a^{z_1} \omega_2(z_2) 
\dots \int _a^{z_{n-1}} \omega_n(z_n).
\eqno(2.8)
$$
These functions of $t$ are holomorphic on $H^-$ and extend holomorphically
to a neighbourhood of  $\bold{P}^1(\bold{R})\setminus \{a,b\}$ (we assume here that $a,b$
are cusps). More precisely, they belong to the $PSL(2,\bold{R})$--module
$\Cal{D}_{\bold{v},-\bold{k}}^{\omega^0,\infty}$ defined in Sec.~1.6 of [BrChDi],
where  $\bold{v}, \bold{k}$ are defined by
$$
\bold{v} (\gamma):= v_1(\gamma)  v_2(\gamma) \dots  v_n(\gamma), \quad
\bold{k} = k_1+k_2 +\dots + k_n.
$$

\smallskip

The key role  in our constructions is played by the following result ([BrCh], Lemma 3.2):

\smallskip

{\bf 2.4.1. Lemma.} {\it The iterated period integral (2.8) as a function of $t \in H$ satisfies for all $\gamma \in SL(2,\Z)$ the following functional equations:
$$
(I_a^b (\omega_1,\dots ,\omega_n; \cdot )|_{-\bold{k}}^{\bold{v}}\gamma )(t)= 
I_{\gamma^{-1}a}^{\gamma^{-1}b} (\omega_1,\dots ,\omega_n; t)
\eqno(2.9)
$$
(see (2.6)).}

\bigskip

\medskip

\centerline{\bf 3. Generalized reciprocity functions from iterated period integrals}

\medskip

{\bf 3.1. Non--commutative generating series.} Fix a finite family of cusp forms as in 2.4
and the respective family of 1--forms $\omega := (\omega_j(z;t))$.
 Let $(A_j),  j=1,\dots ,l$, be independent associative but non--commuting formal variables. 

\smallskip

As in [Ma3], we produce the generating series of all integrals of the type (2.8):
$$
J_a^b(\Omega ;t) := 1+ \sum_{n\ge 1}   \sum_{1\le m_1,\dots ,m_n\le l} I_a^b (\omega_{m_1},\dots ,\omega_{m_n};t)
A_{m_1}\dots A_{m_n}.
$$
\smallskip

Consider the multiplicative group $\bold{G}$ of the formal series in $(A_j)$ with coefficients in functions
of $t$ and lower term 1 ($A_j$ commute with coefficients). The right action of $SL(2,\Z)$ upon this group is defined
coefficientwise. In particular, the action upon $J_a^b(\Omega ;t)$ is given by:
$$
(J_a^b(\Omega ;\cdot)|_{-\bold{k}}^{\bold{v}}\gamma )(t)= := 1+ \sum_{n\ge 1}   \sum_{1\le m_1,\dots ,m_n\le l} (I_a^b (\omega_{m_1},\dots ,\omega_{m_n};\cdot )|_{-\bold{k}}^{\bold{v}}\gamma )(t)A_{m_1}\dots A_{m_n}.
$$

\smallskip

Moreover, the linear action of the group $GL (l,\C )$ upon $(A_1,\dots ,A_l)$ extends
to its action upon the group of formal series $\G$. We will need here only the action of the
subgroup of diagonal matrices.

\smallskip

Let $\gamma \in PSL (2,\bold{Z})$ and $(p,q)\in W.$ Then we will denote by
$\bold{v}(\gamma )*$ the automorphism  of such a ring
sending $A_m$ to  $v_m(\gamma )^{-1}A_m$. 

\smallskip

{\bf 3.1.1. Definition.} {\it The generalised reciprocity function associated to
the family $\omega := (\omega_j)$ is the map of the set
of coprime pairs of positive integers $(p,q)\mapsto f_{\omega}(p,q)\in \bold{G}$
defined by
$$
f_{\omega}(p,q) =f(p,q):=  1+ \sum_{n\ge 1}   \sum_{ 1\le m_1,\dots ,m_n\le l} p^{k_{m_1}+\dots +k_{m_n}}
I_0^{\infty} (\omega_{m_1},\dots ,\omega_{m_n}; qp^{-1})A_{m_1}\dots A_{m_n}.
$$}

Put now
$$
\sigma = \left(\matrix 0& -1\\ 1& 0\endmatrix \right),\quad
\theta = \left(\matrix 1& 1\\ 0& 1\endmatrix \right)\, .
$$

\medskip

{\bf 3.2. Theorem.} {\it The function $f$
satisfies the following functional equation generalising (1.8):
$$
f(p,q) =  (\bold{v}(\theta )* f(p,q+p)) \cdot (\bold{v}(\theta\sigma\theta )* f(q+p,q) ),
\eqno(3.1)
$$
It can be extended to the function defined on coprime pairs $(p,q), p>0, q<0,$
using the functional equation generalising (1.7)
$$
f(p,q):= (\bold{v}(\sigma )*f(-q,p))^{-1}
\eqno(3.2) 
$$
and defining the right hand side via 3.1.1.}

}
\medskip
{\it Remark.}  According to [BrChDi], (2.12), we have
$$
v_m(\sigma )= e^{-\pi  i (k_m+2)/2} ,\     v_m(\theta )= e^{\pi i (k_m+2)/6}.
$$
Hence $v_m(\theta\sigma\theta )=  e^{-\pi i (k_m+2)/6}$.

\smallskip

\medskip

{\bf Proof.} The direct computation shows that
$$
\theta^{-1}(0)=-1,\ \theta^{-1}(\infty )= \infty ,\quad   (\theta\sigma\theta )^{-1}(0)= 0,\
(\theta\sigma\theta )^{-1}(\infty )= -1.
$$
Therefore
$$
J_0^{\infty }(\Omega ;t) = J_{-1}^{\infty }(\Omega ;t) J_0^{-1 }(\Omega ;t) =
J_{\theta^{-1}(0)}^{\theta^{-1}(\infty ) }(\Omega ;t)\cdot J_{(\theta\sigma\theta )^{-1}(0)}^{(\theta\sigma\theta )^{-1}(\infty) }(\Omega ;t) .
\eqno(3.3) 
$$
We now apply (2.9) to both factors in the right hand side of (3.3) term wise. To avoid cumbersome notations, we will
restrict ourselves to simple integrals, that is, to the coefficients of terms linear in $(A_j)$. The general
case easily follows from this one.
\smallskip
 We have then from (2.6):
$$
I_{\theta^{-1}(0)}^{\theta^{-1}(\infty ) }(\omega_j (z;t))= v(\theta )^{-1} j(\theta ,t)^{k_j}\cdot I_0^{\infty} (\omega_j (z; \theta ( t)))=
$$
$$
=v(\theta )^{-1} \cdot I_0^{\infty} (\omega_j (z;  t+1))
\eqno(3.4)
$$
and similarly
$$
I_{(\theta\sigma\theta)^{-1}(0)}^{(\theta\sigma\theta)^{-1}(\infty ) }(\omega_j (z;t))= 
v(\theta\sigma\theta )^{-1} j(\theta\sigma\theta ,t)^{k_j}\cdot I_0^{\infty} (\omega_j (z; \theta ( t)))=
$$
$$
v(\theta\sigma\theta )^{-1}\cdot (t+1)^{k_j}\cdot I_0^{\infty} (\omega_j (z; \frac{t}{t+1})) .
\eqno(3.5)
$$
Thus, from (3.3) we get the following identity for $A_j$--coefficients:
$$
I_0^{\infty}(\omega_j(z;t)) = v(\theta )^{-1}  I_0^{\infty}(\omega_j(z;t+1)) +
v(\theta\sigma\theta )^{-1}\cdot (t+1)^{k_j}\cdot I_0^{\infty} (\omega_j (z; \frac{t}{t+1})).
\eqno(3.6)
$$
In view of growth estimates of period functions in [BrCh] and [BrChDi], we may put here $t=q/p$,
where $q, p$ are coprime integers, $p >  0$. We get

\smallskip

$$
I_0^{\infty}(\omega_j(z;qp^{-1})) = 
$$
$$
v(\theta )^{-1}  I_0^{\infty}(\omega_j(z;(q+p)p^{-1})) +
v(\theta\sigma\theta )^{-1}\cdot ((q+p)p^{-1})^{k_j}\cdot I_0^{\infty} (\omega_j (z; q(q+p)^{-1}))
\eqno(3.7)
$$
that is,
$$
p^{k_j}I_0^{\infty}(\omega_j(z;qp^{-1})) = 
$$
$$
v(\theta )^{-1} p^{k_j} I_0^{\infty}(\omega_j(z;(q+p)p^{-1})) +
v(\theta\sigma\theta )^{-1}\cdot (q+p)^{k_j}\cdot I_0^{\infty} (\omega_j (z; q(q+p)^{-1})).
\eqno(3.8)
$$

Thus, if we put, interpolating  formula (1.4) (established for  even integer weights)
$$
f_{0,j}(p,q):=p^{k_j}I_0^{\infty}(\omega_j(z;qp^{-1})),
\eqno(3.9)
$$
then we get from (3.7) functional equations
$$
f_{0,j}(p,q) =  v(\theta )^{-1} f_{0,j}(p,q+p) + v(\theta\sigma\theta)^{-1} f_{0,j}(q+p,q) .
\eqno(3.10)
$$
This is the $A_j$--linear part of (3.1). The general case is obtained in the same way.
\smallskip
Similarly, looking at the $A_j$--coefficients in the identity $J_0^{\infty }(\Omega ;t)J_{\sigma(0)}^{\sigma ( \infty )}
(\Omega ; t)=1$ and applying (2.6) for $t=-pq^{-1}$ we get
$$
 f_{0,j}(p,q) + v(\sigma) f_{0,j}(-q,p)=0.
$$
The reasoning as above completes the proof.

\smallskip

{\it Remark.} In this statement 3.2, we avoided the simultaneous direct treatment of all cusps $qp^{-1}$
because in our context we have to use the identities of the type $p^k\cdot ((p+q)p^{-1})^k= (p+q)^k$ which
require a separate treatment depending on the signs of integers involved, and the cusps $0$ and $\infty$
should also be treated separately already in the definition of $f(p,q).$ We leave this as an exercise
for the reader.

\medskip

{\bf 3.3. Dedekind cocycles.} In the last section of  [Ma5], equations for  reciprocity functions of even integer weights
(cf. Definition 1.3.1 above) were interpreted as defining a special class of 1--cocycles of $\Gamma := PSL (2,\bold{Z})$.

\smallskip

More precisely, let  $\G$ be a possibly noncommutative left $\Gamma$--module.
It is known that $PSL(2,\bold{Z})$ is the free product of
its two subgroups $\bold{Z}_2$ and $\bold{Z}_3$ generated respectively by
$$
\sigma = \left(\matrix 0& -1\\ 1& 0\endmatrix \right),\quad
\tau = \left(\matrix 0& -1\\ 1& -1\endmatrix \right)\, .
$$

\medskip

Restriction to  $(\sigma ,\tau )$
of any cocycle in  $Z^1(PSL(2,\bold{Z}), \G)$
belongs to the set
$$
\{\, (X,Y) \in \G\times \G\,|\, X\cdot\sigma X=1,\,
Y\cdot\tau Y\cdot\tau^2Y=1\,\}.
\eqno(3.11)
$$

\smallskip

In fact, this defines a bijection between cocycles and pairs (3.11).
\smallskip
An element $(X,Y)$ of (3.11) is called (the representative of)  {\it the Dedekind cocycle}, iff
it satisfies the relation
$$
Y=\tau X. \eqno(3.12)
$$

 Let now $\G_0$ be a group. Denote by $\G$
the group of functions $f:\,   \P^1(\Q)\to \G_0$ with pointwise multiplication. Define the left
action of $\Gamma$ upon $\bold{G}$ by
$$
(\gamma f)(x) = f(\g^{-1}x);\quad f\in \G,\ x\in \P^1(\Q),\ \gamma \in \Gamma.
\eqno(3.13)
$$

Let $f:\, W\to \G_0$ be a $\G_0$--valued reciprocity function, as in Definition 1.3.1.
Define elements $X_f, Y_f\in \G$ as the following functions $\P^1(\Q)\to \G_0$:
$$
X_f(qp^{-1}):= f(p,q),
\eqno(3.14)
$$
$$
Y_f(qp^{-1}):= (\tau X_f)(qp^{-1})= X_f(\tau^{-1}(qp^{-1}))= f(q, q-p).
\eqno(3.15)
$$
\medskip

Then    the map $f\mapsto (X_f,Y_f)$ establishes a bijection between the
set of $\G_0$--valued reciprocity functions and the set of (representatives of)  Dedekind cocycles
from $Z^1(\Gamma ,\G )$  ([Ma5], Theorem 3.6).

\smallskip

We will now show that generalised Dedekind cocycles also can be constructed from
iterated integrals of cusp forms of real weights, although the respective $\Gamma$--module
of coefficients will not be of the form (3.13).

\medskip

{\bf 3.4. A digression: left vs. right.}  Treating cocycles with non commutative coefficients,
we may prefer to work with left or right modules
of coefficients, depending on the concrete environment.  We will give a description of 
Dedekind cocycles with coefficients in a right module.

\smallskip

The  sets of left/right 1--cocycles with coefficients in a left/right module $\G$ are defined by
$$
Z^1_l(\Gamma ,\G) := \{ \lambda :\,\Gamma \to \G\,|\,\lambda (\gamma_1\gamma_2)=\lambda (\g_1)\cdot 
\gamma_1 \lambda(\gamma_2)\,\} .
\eqno(3.16)_l
$$
$$
Z^1_r(\Gamma ,\G) := \{\, \rho :\,\Gamma \to \G\,|\,\rho (\gamma_1\gamma_2)=(\rho (\g_1)|\g_2)\cdot \rho(\g_2)\,\} .
\eqno(3.16)_r
$$
where $*|\gamma$ denotes the right action of $\gamma$.

\smallskip

 Let $\rho \in Z^1_r(PSL(2,\bold{Z}), \G)$ and set $U:=\rho(\sigma ), V:=\rho (\tau )$.
 From  $\sigma^2=\tau^3 =1$ and from  $(3.16)_r$ it follows that

$$
 (U|\sigma)\cdot U=1,\,
(V|\tau^2)\cdot (V|\tau)\cdot V=1.
\eqno(3.17)
$$

\smallskip

Moreover, any pair $(U,V)\in \G\times\G$ satisfying (3.17) comes from a unique cocycle
$\rho \in Z^1_r(PSL(2,\bold{Z}), \G)$.

\smallskip

\medskip

{\bf 3.4.1. Lemma.} {\it (a) If the map  $\G\times\Gamma \to \G$: $(g, \gamma )\mapsto 
g| \gamma $ defines on $\G$ the structure of right $\Gamma$--module, then the map
$( \gamma , g)\mapsto \gamma g:=
g| \gamma^{-1} $ defines on $\G$ the structure of left $\Gamma$--module. 

This construction establishes a bijection between the sets of structures of left/right $\Gamma$--modules on
$\G$.

\smallskip

(b) For such a pair of left/right structures and an element $\lambda \in Z^1_l(PSL(2,\bold{Z}), \G)$,
define  
$$
\rho :\Gamma \to \G :\ \rho (\gamma ) = (\lambda (\gamma^{-1}))^{-1}.
\eqno(3.18)
$$
This establishes a bijection of the respective set of 1--cocycles  $Z^1_l(\Gamma ,\G)$
and $Z^1_r(\Gamma ,\G).$

\smallskip

(c) Assume now that
$\Gamma =PSL(2,\bold{Z})$.  If $(\lambda , \rho )$ is a pair of cocycles connected by (3.18),
put as above $U:=\rho(\sigma ), V:=\rho (\tau )$, and as in [Ma5], sec. 3,  $X:=\lambda(\sigma ), Y:=\lambda (\tau )$.

In this way, Dedekind right cocycles defined by the additional condition $V=U|\sigma$ bijectively correspond
to the Dedekind left cocycles defined in [Ma5] by the condition  $Y=\tau X$.}

\medskip

This can be checked by straightforward computations.

\medskip

{\bf 3.5. Dedekind cocycles from cusp forms of real weights.} In this subsection, we take for $\G$
the subgroup of non--commutative series in $A_i$ (with coefficients in functions on $H^-$) generated by all series
of the form $J_a^b(\Omega;t)$, $a,b\in \bold{P}^1(\Q )$ and fixed $\Omega$
 (cf. 3.1).
\smallskip

The left action of  $\Gamma = PSL(2,\bold{Z})$ upon $\G$ is defined by
$$
\gamma J_a^b(\Omega;t) :=  J_{\gamma a}^{\gamma b}(\Omega;t).
$$

Every cusp $a$ defines the left cocycle $\lambda_a: \,\Gamma\to \G$: 
$$
\lambda_a(\gamma ):= J^a_{\gamma a} (\Omega ;t)
$$
because
$$
J^a_{\gamma_1\gamma_2 a} = J^a_{\gamma_1 a} \cdot  J^{\gamma_1 a}_{\gamma_1\gamma_2 a}=
J^a_{ \gamma_1 a} \cdot  \gamma_1 J^a_{\gamma_2 a}.
$$
\medskip

In order to construct our Dedekind cocycle, we combine the $\sigma$--component
of  $\lambda_{\infty}$ with $\tau$--component of $\lambda_0$:

{\bf 3.5.1. Theorem.} {\it The pair
$$
 X:=J_0^{\infty}(\Omega ;t), \quad  Y:=J_1^0(\Omega ;t)
  \eqno(3.19)
 $$ 
 is (the representative of) a left Dedekind cocycle. }
\medskip

{\bf Proof.}  We have
$$
X\cdot \sigma X=  J_0^{\infty}(\Omega ;t)\cdot J^0_{\infty}(\Omega ;t) =1,
$$
$$
Y\cdot\tau Y\cdot\tau^2Y= J_1^0(\Omega ;t)\cdot J^1_{\infty}(\Omega ;t)
\cdot J_0^{\infty}(\Omega ;t) =1,
$$
and finally
$$
\tau X=  J_{\tau 0}^{\tau \infty}(\Omega ;t) = J_1^0(\Omega ;t) =Y.
$$

\newpage

\bigskip
\centerline{\bf References}

\medskip

[Ap] T.~M.~Apostol. {\it Generalized Dedekind sums and transformation formulae of certain Lambert series.} Duke 
Math.~J. 17 (1950), 147--157.

\smallskip

[BaSch]  S.~Baumard, L.~Schneps. {\it On the derivation representation
of the fundamental Lie algebra of mixed elliptic motives.} arXiv:1510.05549

\smallskip

[BrChDi] R.~Bruggeman, Y.~Choie, N.~Diamantis. {\it Holomorphic automorphic
forms and cohomology.} arXiv:1404.6718 .

\smallskip

[BrCh]  R.~Bruggeman, Y.~Choie. {\it Multiple period integrals and cohomology.}
arXiv:1512.06225 .
\smallskip

[ChZ] Y.~J.~Choie, D.~Zagier. {\it Rational period functions
for $PSL(2,\bold{Z} ).$} In: A Tribute to Emil Grosswald:
Number Theory and relatied Analysis, Cont. Math.,
143 (1993), AMS, Providence, 89--108.

\smallskip

[Fu1] Sh.~Fukuhara. {\it Modular forms, generalized Dedekind symbols
and period polynomials.} Math. Ann., 310 (1998), 83--101.

\smallskip

[Fu2] Sh.~Fukuhara. {\it Dedekind symbols with polynomial reciprocity laws.}
Math. Ann., 329 (2004), No. 2,  315--334.

\smallskip

[Hai] R.~Hain. {\it The Hodge--de Rham theory of modular groups.} arXiv:1403.6443

\smallskip

[HaiMo] R.~Hain, M.~Matsumoto. {\it Universal mixed elliptic motives.}  arXiv:1512.03975

\smallskip

[KiMel] R.~C.~Kirby,  P.~Melvin. {\it Dedekind sums, $\mu$--invariants and the signature cocycle.} 
Math. Ann. 299 (1994), 231--267.

\smallskip


[Kn] M.~Knopp. {\it Some new results on the Eichler cohomology
of automorphic forms.} Bull.~AMS, Vol.~80, Nr.~4 (1976), 607--632.

\smallskip


[KnMa]   M.~Knopp, H.~Mawi.  {\it Eichler cohomology theorem for automorphic forms of small weights.}
Proc.~AMS, Vol.~138, Nr.~2 (2010), 395--404.
\smallskip

[Ma1] Yu.~Manin. {\it Parabolic points and zeta-functions of modular curves.}
Russian: Izv. AN SSSR, ser. mat. 36:1 (1972), 19--66. English:
Math. USSR Izvestiya, publ. by AMS, vol. 6, No. 1 (1972), 19--64,
and Selected Papers, World Scientific, 1996, 202--247.

\smallskip

[Ma2] Yu.~Manin. {\it Periods of parabolic forms and $p$--adic Hecke series.}
Russian: Mat. Sbornik, 92:3 (1973), 378--401. English:
Math. USSR Sbornik, 21:3 (1973), 371--393,
and Selected Papers, World Scientific, 1996, 268--290.

\smallskip

[Ma3] Yu.~Manin. {\it Iterated integrals of modular forms and 
noncommutative modular symbols.}
In: Algebraic Geometry and Number Theory.
In honor of V. Drinfeld's 50th birthday.
Ed. V.~Ginzburg. Progress in Math., vol. 253. 
Birkh\"auser, Boston, pp. 565--597. Preprint math.NT/0502576.

\smallskip

[Ma4] Yu.~Manin. {\it Iterated Shimura integrals.} Moscow Math. Journal, 
vol. 5, Nr. 4 (2005), 869--881.
 Preprint math.AG/0507438.

\smallskip

[Ma5]  Yu.~Manin. {\it  Non--commutative generalized Dedekind symbols.} Pure and Appl. Math.
Quarterly, Vol.10, Nr. 2 (2014), 245 -- 258.
 Preprint arXiv:1301.0078
 
 \smallskip
 
 [Ma6]  Yu.~Manin. {\it Local zeta factors and geometries under $\roman{Spec}\,\bold{Z}$.} 9 pp.
Preprint arXiv:1407.4969   (To be published in Izvestiya RAN, issue dedicated to J.-P.~Serre).

\smallskip
[Mel] P.~Melvin. {\it Tori in the diffeomorphism groups of of simply--connected  4--manifolds.}
Math. Proc. Cambridge Phil. Soc. 91 (1982),  305--314.

\smallskip

[Po] A.~Pollack. {\it Relations between derivations arising from modular forms.}
http://dukespace.lib.duke.edu/dspace/handle/10161/1281

\bigskip

\enddocument